\documentstyle[12pt]{article}

\fussy
\begin{document}

\begin{center}
{\large {\bf {Twist deformations for generalized Heisenberg algebras}}}
\vskip
5mm {\sc Vladimir D. Lyakhovsky
\footnote{E-mail address: lyakhovs@snoopy.phys.spbu.ru } }\\ \vskip 2mm
Theoretical Department, St. Petersburg State University, 198904, St.
Petersburg, Russia \\ \vskip 2mm
\end{center}

\begin{abstract}
Multidimensional Heisenberg algebras, whose creation  $a^+$ and
annihilation $a^-$ operators are the n-dimensional vectors, can be injected
into simple Lie algebras $g$. It is demonstrated that the spectrum of
their deformations can be investigated using chains of extended Jordanian
twists applied to $U(g)$. In the case of $U(sl(N))$ (for $N>5$ ) the
two-dimensional Heisenberg subalgebras $\widetilde{\cal H}$ have nine deformed
costructures connected by four "internal" and "external" twists composing the
commutative diagram.

\end{abstract}


\section{Introduction}

A Hopf algebra ${\cal A}(m,\Delta ,\epsilon,S)$ can be transformed
\cite{D2} by an
invertible {\bf twisting element} ${\cal F}\in {\cal A} \otimes {\cal A}$,
into a {\bf twisted algebra} ${\cal A}_{{\cal F}}(m,\Delta _{{\cal F}},
\epsilon ,S_{{\cal F}})$, that has the same multiplication and counit but
the twisted coproduct and the antipode given by
\begin{equation}
\label{def-t} 
\Delta _{{\cal F}}(a)={\cal F}\Delta (a){\cal F}^{-1},
S_{{\cal F}}(a)=vS(a)v^{-1}, v=\sum f_i^{(1)}S(f_i^{(2)}), a\in {\cal
A}.
\end{equation}
The twisting element has to satisfy the equations
\begin{equation}
\label{TE} {\cal F}_{12}(\Delta \otimes id)({\cal F}) = {\cal F}_{23}(id
\otimes \Delta)({\cal F}), \quad (\epsilon \otimes id)({\cal F}) = (id
\otimes \epsilon)({\cal F})=1.
\end{equation}

Let ${\cal A}$ and ${\cal B}$ be the universal enveloping algebras: ${\cal A}
= U(l) \subset {\cal B}= U(g)$ with $l \subset
g$. If $U(l)$ is the minimal subalgebra on which ${\cal F
}$ is completely defined as ${\cal F} \in U(l) \otimes U(l)$ then
$l$ is called the {\bf carrier} algebra for ${\cal F}$.

The first and well known example of the twist that was written in the
explicit form \cite{OGIEV} corresponds to the carrier subalgebra $B(2)$ with
generators $H$ and $E$, $[H,E]=E$. It is called the {\bf Jordanian twist}
and has the twisting element
\begin{equation}
\label{ogiev}\Phi _{\cal J}=e^{H\otimes \sigma },\quad
\quad \sigma =\ln (1+E).
\end{equation}

For the {\bf extended Jordanian twists} suggested in \cite
{KLM} the carrier subalgebra ${\bf L}$ is four-dimensional:
\begin{equation}
\label{l-def}
\begin{array}{l}
[H,E]=E,\quad [H,A]=\alpha A,\quad [H,B]=\beta B, \\
[0.2cm][A,B]=E,\quad [E,A]=[E,B]=0,\quad \quad \alpha +\beta =1.
\end{array}
\end{equation}
The corresponding twisting element contains the {\bf Jordanian factor}:
\begin{equation}
\label{t-ext}
{\cal F}_{{\cal E(\alpha ,\beta )}}=\Phi _{{\cal E(\alpha
,\beta )}}\Phi _{\cal J}
\end{equation}
and the {\bf extension}
\begin{equation}
\label{fractions}
\Phi _{{\cal E(\alpha ,\beta )}}=\exp \{A\otimes Be^{-\beta
\sigma }\}.
\end{equation}
This twist defines the deformed Hopf algebras ${\bf L}_{{\cal E(\alpha
,\beta )}}$ with the costructure
\begin{equation}
\label{e-costr}
\begin{array}{lcl}
\Delta _{{\cal E(\alpha ,\beta )}}\,(H) & = & H\otimes e^{-\sigma }+1\otimes
H-A\otimes Be^{-(\beta +1)\sigma }, \\
\Delta _{{\cal E(\alpha ,\beta )}}\,(A) & = & A\otimes e^{-\beta \sigma
}+1\otimes A, \\
\Delta _{{\cal E(\alpha ,\beta )}}\,(B) & = & B\otimes e^{\beta \sigma
}+e^\sigma \otimes B, \\
\Delta _{{\cal E(\alpha ,\beta )}}\,(E) & = & E\otimes e^\sigma +1\otimes E.
\end{array}
\end{equation}

In general the composition of two twists is not a twist. But there are some
important examples of the opposite behavior. When ${\bf L}$ is a subalgebra
${\bf L\subset }{\mbox{g }}$ there may exist several pairs of generators of
the type $(A,B)$ arranged so that the Jordanian twist can acquire several
similar extensions \cite{KLM}. This demonstrates that some twistings can be
applied successively to the initial Hopf algebra even in the case when their
carrier subalgebras are nontrivially linked. In the universal enveloping
algebras for classical Lie algebras there exists the possibility to
construct systematically the special sequences of twists called {\bf chains}
\cite{KLO}:
\begin{equation}
\label{mainchain}
{\cal F}_{{\cal B}_{p\prec 0}}\equiv {\cal F}_{{\cal B}_p}
{\cal F}_{{\cal B}_{p-1}}\ldots {\cal F}_{{\cal B}_0}.
\end{equation}
The factors ${\cal F}_{{\cal B}_k}=\Phi _{{\cal E}_k}\Phi _{{\cal J}_k}$ of
the chain are the twisting elements of the extended Jordanian twists for the
initial Hopf algebra ${\cal A}_0$. Here the extensions $\left\{ \Phi _{{\cal
E}_k},k=0,\ldots p-1\right\} $ contain the fixed set of normalized factors $
\Phi _{{\cal E(\alpha ,\beta )}}=\exp \{A\otimes Be^{-\beta \sigma }\}$ ,
the so called {\bf full set}. It was proved that in the classical Lie
algebras that conserve symmetric invariant forms such chains can be made
maximal and proper. This means that for the algebras $U(A_n)$, $U(B_n)$ and $
U(D_n)$ there exist chains ${\cal F}_{{\cal B}_{p\prec 0}}$ that cannot be
reduced to a chain for a simple subalgebra and their full sets of
extensions are the maximal sets in the sense described below.

To construct a maximal proper chain for a classical Lie algebra the
sequence ${\cal A}\equiv {\cal A}_0\supset {\cal A}_1\supset
\ldots \supset {\cal A}_{p-1}\supset {\cal A}_p$ of Hopf subalgebras is to
be fixed. For example in the case  of $U(sl(N))$ the corresponding sequence
is:
$
U(sl(N))\supset U(sl(N-2))\supset \ldots \supset U(sl(N-2k))\ldots
$.
In each element of the sequence the {\bf initial root} $\lambda _0^k$ must
be chosen. Here $\lambda _0^k=e_1-e_2$ for each $sl(M)$
(the roots are written in the standard $e$-basis).
For the initial root let us form the set $\pi _k$ of {\bf
constituent roots},
\begin{equation}
\label{kpi}
\pi _k=\left\{ \lambda ^{\prime },\lambda ^{\prime \prime
}|\lambda ^{\prime }+\lambda ^{\prime \prime }=\lambda _0^k;\quad \lambda
^{\prime }+\lambda _0^k,\lambda ^{\prime \prime }+\lambda _0^k \neg \in
\Lambda _{{\cal A}}\right\}
\end{equation}
(here $\Lambda _{{\cal A}}$ is the root system of ${\cal A}_0$). For each
element $\lambda ^{\prime }\in \pi _k$ one can choose such an element $
\lambda ^{\prime \prime }\in \pi _k$ that $\lambda ^{\prime }+\lambda
^{\prime \prime }=\lambda _0^k$. So, $\pi _k$ is naturally decomposed as
\begin{equation}
\pi _k=\pi _k^{\prime }\,\cup \,\pi _k^{\prime \prime },\quad \quad \pi
_k^{\prime }=\{\lambda ^{\prime }\},\quad \pi _k^{\prime \prime }=\{\lambda
^{\prime \prime }\}.
\end{equation}

In these terms the factors ${\cal F}_{{\cal B}_k}$ of the chain (\ref
{mainchain}) are fixed as follows:
\begin{equation}
\label{fabk}{\cal F}_{{\cal B}_k}=\Phi _{{\cal E}_k}\Phi _{{\cal J}_k}
\end{equation}
with
\begin{equation}
\Phi _{{\cal J}_k}=\exp \{H_{\lambda _0^k}\otimes \sigma _0^k\},\quad \quad
\sigma _0^k=\ln (1+L_{{\scriptsize {\lambda _0^k}}});
\end{equation}
\begin{equation}
\Phi _{{\cal E}_k}=\prod_{\lambda ^{\prime }\in \pi _k^{\prime }}\Phi _{
{\cal E}_{\lambda ^{\prime }}}=\prod_{\lambda ^{\prime }\in \pi _k^{\prime
}}\exp \{L_{\lambda ^{\prime }}\otimes L_{{\scriptsize {\lambda _0^k}
-\lambda ^{\prime }}}e^{-\frac 12\sigma _0^k}\}
\end{equation}
(here $L_\lambda $ is the generator associated to the root $\lambda $).

\section{Deformed Heisenberg algebras}

One of the characteristic features of the extension $\Phi _{{\cal E}}$ in
the ordinary extended Jordanian twists is that it connects the
Heisenberg subalgebras ${\cal H}_{{\cal J}}$ and ${\cal H}_{{\cal EJ}}$ with
$[A,B]=E, \alpha +\beta =1 $ and
different deformed costructures (see \cite{VM}):
\begin{equation}
\label{pureheis}
\begin{array}{ll}
\Phi _{{\cal E}}:\left\{
\begin{array}{lcl}
\Delta _{{\cal J}}(A) & = & A\otimes e^{\alpha \sigma }+1\otimes A, \\
\Delta _{{\cal J}}(B) & = & B\otimes e^{\beta \sigma }+1\otimes B, \\
\Delta _{{\cal J}}(E) & = & E\otimes e^\sigma +1\otimes E,
\end{array}
\right\} & \longrightarrow \\
 \left\{
\begin{array}{lcl}
\Delta _{{\cal EJ}}\,(A) & = & A\otimes e^{-\beta \sigma }+1\otimes A, \\
\Delta _{{\cal EJ}}\,(B) & = & B\otimes e^{\beta \sigma }+e^\sigma \otimes
B,  \\
\Delta _{{\cal EJ}}(E) & = & E\otimes e^\sigma +1\otimes E.
\end{array}
\right\} &
\end{array}
\end{equation}

A chain of twists (\ref{mainchain}) ${\cal F}_{{\cal B}_{p\prec 0}}
\equiv {\cal F}_{{\cal B}_p}
{\cal F}_{{\cal B}_{p-1}}\ldots {\cal F}_{{\cal B}_0}
:{\cal A} \longrightarrow {\cal A}_{{\cal B}_{p\prec 0}}$
contains $p + 1$ Jordanian factors $\Phi _{{\cal J}_k}$. Their product $\Phi
_{{\cal J}_{p\prec 0}}\equiv \prod_k{\Phi }_{{\cal J}_{k-1}}$ is also a
twisting element for the universal enveloping algebra ${\cal A}$,
it produces the {\bf multijordanian deformation}
$
\Phi _{{\cal J}_{p\prec 0}}:
{\cal A}\longrightarrow {\cal A}_{{\cal J}_{p\prec 0}}
$.
This means that the product $\widehat{{\Phi }}_{{\cal E}_{p\prec 0}}
\equiv \prod_k \widehat{{\Phi }}_{{\cal E}_k}$ with
$
\widehat{{\Phi }}_{{\cal E}_k}=(\prod_{m>k}{\Phi }_{{\cal J}_m}){\Phi }_{
{\cal E}_k}(\prod_{m>k}{\Phi }_{{\cal J}_m})^{-1}
$
is in turn a twisting element for ${\cal A}_{{\cal J}_{p\prec 0}},$
$
\widehat{\Phi }_{{\cal E}_{p\prec 0}}:{\cal A}_{{\cal J}_{p\prec
0}}\longrightarrow {\cal A}_{B_{p\prec 0}}
$.

In \cite{LO-3} it was shown that $\widehat{\Phi }_{{\cal E}_{p\prec 0}}$
plays the role of an extension for the multijordanian twists $\Phi _{{\cal J}
_{p\prec 0}}$ . It was proved that there are the subalgebras in ${\cal A}$
that are carriers for $\widehat{\Phi }_{{\cal E}_{p\prec 0}}$ and whose
costructures are shifted  from one possible deformed ''state'' to the
other by the twists $\widehat{\Phi }_{{\cal E}_{p\prec 0}}$.

Here we consider a certain type of such subalgebras, namely the
multidimensional Heisenberg subalgebras $\widetilde{\cal H}(2,N-4)$
in $U\left(sl(N)\right) $ generated by the $2\times 2$-block of
central and $2\times
(N-4)$ pairs of creation and annihilation operators:
\begin{equation}
\label{generators}
\begin{array}{ll}
\frame{$
\begin{array}{cccc}
E_{13} & E_{14} & \ldots  & E_{1,N-2} \\
E_{23} & E_{24} & \ldots  & E_{2,N-2}
\end{array}
$} & \frame{$
\begin{array}{cc}
E_{1,N-1} & \quad E_{1N}\quad  \\
E_{2,N-1} & \quad E_{2N}\quad
\end{array}
$} \\  & \frame{$
\begin{array}{cc}
E_{3,N-1} & E_{3N} \\
\vdots  & \vdots  \\
E_{N-3,N-1} & E_{N-3,N} \\
E_{N-2,N-1} & E_{N-2,N}
\end{array}
$}
\end{array}
\end{equation}
We shall normalize the Cartan elements as $H_{i,k}=1/2(E_{ii}-E_{kk})$, use
the standard $gl(N)$-basis $\{E_{ij}\}_{i,j=1,\dots N}$, and $\sigma =\ln
(1+E)$.

First we shall study the properties
of $\widetilde{\cal H}_{\cal F}(2,N-4)$ twisted by the factors of the
2-chain:
\begin{equation}
\label{2-chain}{\cal F}_{{\cal B}_{1\prec 0}}
={\Phi }_{{\cal E}_1}{\Phi }_{
{\cal J}_1}{\Phi }_{{\cal E}_0}{\Phi }_{{\cal J}_0}
\end{equation}
with
\begin{equation}
\begin{array}{lcl}
{\Phi }_{{\cal J}_{k-1}}&=& \exp (H_{k,N-k+1}\otimes \sigma _{k,N-k+1}), \\ {
\Phi }_{{\cal E}_{k-1}}&=&\exp \left( \sum_{s=k+1}^{N-k}E_{k,s}\otimes
E_{s,N-k+1}e^{-\frac 12\sigma _{k,N-k+1}}\right) \\
&=& \prod_{r=k+1}^{N-k}{\Phi }_{{\cal E}_{k-1}(r)} .
\end{array}
\end{equation}
To visualize the resulting deformations the following notations will be used
for the costructures:
\begin{equation}
\label{coblock}
\begin{array}{c}
\begin{array}{lll}
P^{0}\left( L\right) & \equiv & L \otimes 1 + 1 \otimes L, \\
P_i^{\pm }\left( L\right) & \equiv & L\otimes e^{\pm \frac 12\sigma
_{i,N+1-i}}+1\otimes L, \\
R_i\left( L\right) & \equiv & L\otimes e^{\frac 12\sigma
_{i,N+1-i}}+e^{\sigma _{i,N+1-i}}\otimes L, \\
T_i\left( L\right) & \equiv & L\otimes e^{\sigma _{i,N+1-i}}+1\otimes L, \\
T^{++}\left( L\right) & \equiv & L\otimes e^{\frac 12\sigma _{1,N}+\frac
12\sigma _{2,N-1}}+1\otimes L, \\
T^{\mp \pm }\left( L\right) & \equiv & L\otimes e^{\mp \frac 12\sigma
_{1,N}\pm \frac 12\sigma _{2,N-1}}+1\otimes L, \\
T_{R_i}\left( L\right) & \equiv & L\otimes e^{\frac 12\sigma _{1,N}+\frac
12\sigma _{2,N-1}}+e^{\sigma _{i,N+1-i}}\otimes L, \\
S_{2}^{-} & \equiv & -E_{2,r}\otimes E_{1,N-1}e^{-\frac 12\sigma _{2,N-1}},
\\
S_{2}^{+} & \equiv & E_{2,N}\otimes E_{r,N-1}e^{\frac 12\sigma _{1,N}-\frac
12\sigma _{2,N-1}}, \\
S_{1}^{-} & \equiv & -E_{1,r}\otimes E_{2,N}e^{-\frac 12\sigma _{1,N}}, \\
S_{1}^{+} & \equiv & E_{1,N-1}\otimes E_{r,N}e^{-\frac 12\sigma _{1,N}+\frac
12\sigma _{2,N-1}},
\end{array}
\\
i=1,2;\quad r=3,\ldots ,N-2.
\end{array}
\end{equation}
In these terms the transformations performed in ${\cal H}$ by the factors of
the canonical extended Jordanian twist (use the formulas
(\ref{l-def})-(\ref{e-costr}) with $\alpha = \beta = 1/2$ ) can be written
as
\begin{equation}
\label{pureheisjord}
\begin{array}{c}
\Phi _{{\cal J}}:\left\{
\begin{array}{lc}
P^0(A) & P^0 \left( E \right) \\
& P^0 (B)
\end{array}
\right\} \longrightarrow \left\{
\begin{array}{lc}
P^+ (A) & T \left( E\right) \\
& P^+(B)
\end{array}
\right\}.
\end{array}
\end{equation}

\begin{equation}
\label{pureheisco}
\begin{array}{c}
\Phi _{{\cal E}}:\left\{
\begin{array}{lc}
P^{+}(A) & T \left( E \right) \\
& P^{+} (B)
\end{array}
\right\} \longrightarrow \left\{
\begin{array}{lc}
P^{-} (A) & T \left( E\right) \\
& R(B)
\end{array}
\right\}.
\end{array}
\end{equation}

Let us perform the 2-Jordanian twisting in $\widetilde{\cal H}(2,N-4)$,
$
\Phi _{{\cal J}_1{\cal J}_2}:\widetilde{\cal H}(2,N-4)
\longrightarrow \widetilde{\cal H}_{{\cal J}_1{\cal J}_2}(2,N-4)
$.
Here the deformed coproducts will be presented by the schemes analogous
to (\ref{pureheisjord}) and (\ref{pureheisco}). In terms
of (\ref{coblock}) the costructure
of $\widetilde{\cal H}_{{\cal J}_1{\cal J}_2}(2,N-4)$ will acquire
the following form:
\begin{equation}
\label{2jordco}
\begin{array}{ll}
\frame{$
\begin{array}{ccc}
P_1^{+}\left( E_{13}\right)  &  \ldots  &
P_1^{+}\left( E_{1,N-2}\right)  \\
P_2^{+}\left( E_{23}\right)  &  \ldots  &
P_2^{+}\left( E_{2,N-2}\right)
\end{array}
$} & \frame{$
\begin{array}{cc}
T^{++}\left( E_{1,N-1}\right)  & \quad T_1\left( E_{1N}\right)  \\
T_2\left( E_{2,N-1}\right)  & \quad T^{++}\left( E_{2N}\right)
\end{array}
$} \\  & \frame{$
\begin{array}{cc}
P_2^{+}\left( E_{3,N-1}\right)  & P_1^{+}\left( E_{3N}\right)  \\
\vdots  & \vdots  \\
P_2^{+}\left( E_{N-2,N-1}\right)  & P_1^{+}\left( E_{N-2,N}\right)
\end{array}
$}
\end{array}
\end{equation}
As we have mentioned above any Heisenberg subalgebra with the costructure $
\left\{ P^{+},T,P^{+}\right\} $ can be twisted by the corresponding
extension $\Phi _{{\cal E}}$ . In (\ref{2jordco}) such are the triples $
\left\{ P_i^{+}\left( E_{ir}\right) \right. $, $ T_i\left( E_{i,N-i+1}\right)
$, $\left. P_i^{+}\left( E_{r,N-i+1}\right) \right\}$. Obviously the
deformations performed by the extensions $\Phi _{{\cal E}_{i-1}(r)}
=\exp \left(E_{i,r}\otimes E_{r,N-i+1}e^{-\frac 12\sigma _{i,N-i+1}}\right)$
in these triples do not touch the generators other than $\left\{
E_{1,r},E_{2,r}\right. , \left. E_{r,N-1},E_{r,N}\right\} $.
Thus to study the deformations induced on $ \widetilde{\cal H}(2,N-4)$ by
the chain of twists (\ref{2-chain}) it is sufficient to focus the attention
on the behavior of the subalgebra $\widetilde{\cal H} \equiv
\widetilde{\cal H}(2,1),$
\begin{equation}
\label{rgenerators}
\begin{array}{cc}
\frame{$
\begin{array}{c}
E_{1r} \\
E_{2r}
\end{array}
$} & \frame{$
\begin{array}{cc}
E_{1,N-1} & E_{1N} \\
E_{2,N-1} & E_{2N}
\end{array}
$} \\  & \frame{$
\begin{array}{cc}
E_{r,N-1} & E_{rN}
\end{array}
$}
\end{array}
\end{equation}
Starting with its 2-Jordanian deformation
$ \widetilde{\cal H}_{{\cal J}_1{\cal J}_0}$ one can apply any of the two
extensions: $\Phi _{{\cal E}_{i-1}(r)}=\exp \left( E_{i,r}
\otimes E_{r,N-i+1}e^{-\frac12\sigma _{i,N-i+1}}\right) $
$\left( i=1,2\right) $. Notice that in the chain (\ref{2-chain}) the
extension $\Phi _{{\cal E}_1(r)}$ commutes not only with all the other
extension factors $\Phi _{{\cal E}_{i-1}(r)}$ but
also with ${\Phi }_{{\cal J}_1}$ .We get three twisted ''states'':
$$
\widetilde{\cal H}_{{\cal E}_0{\cal J}_1{\cal J}_0}
\stackrel{\Phi _{{\cal E}_0(r)}}{\longleftarrow } \widetilde{\cal H}
_{{\cal J}_1{\cal J}_0}\stackrel{\Phi _{{\cal E}_1(r)}}{
\longrightarrow } \widetilde{\cal H}_{{\cal E}_1{\cal J}_1{\cal J}_0}
$$
The corresponding coalgebras are transformed as follows
\begin{equation}
\label{r2exts}
\begin{array}{c}
\left\{
\begin{array}{ll}
\frame{$
\begin{array}{c}
P_1^{-} \\
P_2^{+}+S_1^{-}
\end{array}
$} & \frame{$
\begin{array}{cc}
T^{++} & \quad T_1\quad  \\
T_2 & T^{++}
\end{array}
$} \\  & \frame{$
\begin{array}{cc}
P_2^{+}+S_1^{+} & R_1
\end{array}
$}
\end{array}
\right\}
\longleftarrow
\left\{
\begin{array}{ll}
\frame{$
\begin{array}{c}
P_1^{+} \\
P_2^{+}
\end{array}
$} & \frame{$
\begin{array}{cc}
T^{++} &  T_1 \\
T_2 & T^{++}
\end{array}
$} \\  & \frame{$
\begin{array}{cc}
P_2^{+} & \;\;P_1^{+}
\end{array}
$}
\end{array}
\right\}
\\
\longrightarrow
\left\{
\begin{array}{ll}
\frame{$
\begin{array}{c}
P_1^{+}+S_2^{-} \\
P_2^{-}
\end{array}
$} & \frame{$
\begin{array}{cc}
T^{++} & \quad T_1 \\
T_2 & \quad T^{++}
\end{array}
$} \\  & \frame{$
\begin{array}{cc}
R_2 & P_1^{+}+S_2^{+}
\end{array}
$}
\end{array}
\right\}
\end{array}
\end{equation}

Now we shall show that among the twists for $U\left( sl(N)\right) $ one can
find those that perform further deformations of the states (\ref{r2exts}).
Their carrier subalgebras are not included in the $\widetilde{\cal H}$, so
with respect to $\widetilde{\cal H}$ these twists must be treated as
external. Consider again the chain (\ref{2-chain}). According to the
''matreshka'' effect (see \cite{KLO}) after the first two twists ${\Phi }_{
{\cal J}_0}$ and ${\Phi }_{{\cal E}_0}$ the generators in the subalgebra $U_{
{\cal E}_0{\cal J}_0}\left( sl(N-2)\right) $ will acquire trivial coproducts.
Consequently when the first three factors ${\Phi }_{{\cal J}_1}{\Phi }_{
{\cal E}_0}{\Phi }_{{\cal J}_0}$ are applied the generators $\left\{
E_{2,r},E_{r,N-1}\right\} $ return to the state $P_2^{+}$. The twisting
factor that performs this transition appears when we drag the Jordanian
twisting element ${\Phi }_{{\cal J}_1}$ to the right:
\begin{equation}
\label{adjoint0}
\begin{array}{c}
{\Phi }_{{\cal J}_1}{\Phi }_{{\cal E}_0}{\Phi }_{{\cal J}_0} ={\Phi }_{{\cal
J}_1}{\Phi }_{{\cal E}_0(2)}{\Phi }_{{\cal E}_0(N-1)}{\Phi }_{{\cal E}
_0(N-2\prec 3)} {\Phi }_{{\cal J}_0}=
\\
{\Phi }_{{\cal J}_1}{\Phi }_{{\cal E
}_0(2)}{\Phi }_{{\cal E}_0(N-1)} {\Phi }_{{\cal J}_1}^{-1}{\Phi }_{{\cal J}
_1} {\Phi }_{{\cal E}_0(N-2\prec 3)}{\Phi }_{{\cal J}_0}=
\\
\widetilde{\Phi }_{
{\cal E}_0(2,N-1)} {\Phi }_{{\cal E}_0(N-2\prec 3)}{\Phi }_{{\cal J}_1{\cal J
}_0.}
\end{array}
\end{equation}
The factor ${\Phi }_{{\cal E}_0(N-2\prec 3)}$ is a twisting element for $U_{
{\cal J}_1{\cal J}_0}$. The relation (\ref{adjoint1}) signifies that
\begin{equation}
\label{adjoint1}
\begin{array}{c}
\widetilde{\Phi }_{{\cal E}_0(2,N-1)}=\exp \left( \left( E_{1,2}+\frac
12E_{1,N-1}+E_{1,N-1}H_{2,N-1}\right)\otimes
\right. \\ \left.
\otimes E_{2,N}e^{-\frac 12\left(
\sigma _{1,N}+\sigma _{2,N-1}\right) } +E_{1,N-1}\otimes
E_{N-1,N}e^{-\frac 12\left( \sigma _{1,N}-\sigma _{2,N-1}\right) }\right)
\end{array}
\end{equation}
twists the costructure of $U_{{\cal E}_0(N-2\prec 3){\cal J}_1{\cal J}_0}$.
It is important that for the generators in (\ref{adjoint1}) the coproducts
in $U_{{\cal J}_1{\cal J}_0}$ and in $U_{{\cal E}_0(N-2\prec 3){\cal J}_1
{\cal J}_0}$ are the same. This means that the factor $\widetilde{\Phi }_{
{\cal E}_0(2,N-1)}$ can twist directly the algebra $U_{{\cal J}_1{\cal J}_0}$
and the subalgebra $\widetilde{\cal H}_{{\cal J}_1{\cal J}
_0}$ in it. The structure constants of this subalgebra are invariant with
respect to the renumbering $\left( 1\rightleftharpoons
2,N-1\rightleftharpoons N\right) $ of the indices of generators. Thus there
exists the second external twisting factor
\begin{equation}
\label{adjoint2}
\begin{array}{c}
\widetilde{\Phi }_{{\cal E}_1(2,N-1)}=\exp \left( \left( E_{2,1}+\frac
12E_{2,N}+E_{2,N}H_{1,N}\right) \otimes
\right. \\ \left.
\otimes E_{1,N-1}e^{-\frac 12\left( \sigma
_{1,N}+\sigma _{2,N-1}\right) }+E_{2,N}\otimes
E_{N,N-1}e^{+\frac 12\left( \sigma _{1,N}-\sigma _{2,N-1}\right) }\right)
\end{array}
\end{equation}
that can be applied to $\widetilde{\cal H}_{{\cal J}_1
{\cal J}_0}$. One can find this twisting element directly considering the
chain
\begin{equation}
\label{addchain}
{\cal F^{\prime }}_{{\cal B}_{0\prec 1}}={\Phi ^{\prime }}_{
{\cal E}_0}{\Phi }_{{\cal J}_0}
{\Phi ^{\prime }}_{{\cal E}_1}{\Phi }_{{\cal J}_1}
\end{equation}
where the maximal set of the constituent roots is used in the extension ${
\Phi ^{\prime }}_{{\cal E}_1}$ (that is for the root $\left(
e_2-e_{N-1}\right) $). The algebra $U_{{\cal J}_1{\cal J}_0}= U_{{\cal J}_0
{\cal J}_1}$ can be considered as an intermediate deformed object for both
chains (\ref{2-chain}) and (\ref{addchain}).

When any of the external factors $\widetilde{\Phi }_{{\cal E}_0(2,N-1)}$
and $\widetilde{\Phi }_{{\cal E}_1(2,N-1)}$ are applied
to $\widetilde{\cal H}_{{\cal J}_1{\cal J}_0}$ one of the pairs of
creation-annihilation generators retains the costructure $P^{+}$:
\begin{equation}
\label{r2extexts}
\begin{array}{c}
\left\{
\begin{array}{ll}
\frame{$
\begin{array}{c}
P_1^{+} \\
P_2^{+}-S_1^{-}
\end{array}
$} & \frame{$
\begin{array}{cc}
T^{-+} & \quad T_1\quad  \\
T_2 & T_{R_1}
\end{array}
$} \\  & \frame{$
\begin{array}{cc}
P_2^{+}-S_1^{+} & P_1^{+}
\end{array}
$}
\end{array}
\right\}
\leftarrow
\left\{
\begin{array}{ll}
\frame{$
\begin{array}{c}
P_1^{+} \\
P_2^{+}
\end{array}
$} & \frame{$
\begin{array}{cc}
T^{++} & T_1 \\
T_2 & T^{++}
\end{array}
$} \\  & \frame{$
\begin{array}{cc}
P_2^{+} & \;\;P_1^{+}
\end{array}
$}
\end{array}
\right\}
\\
\longrightarrow
\left\{
\begin{array}{ll}
\frame{$
\begin{array}{c}
P_1^{+}-S_2^{-} \\
P_2^{+}
\end{array}
$} & \frame{$
\begin{array}{cc}
T_{R_2} & \quad T_1 \\
T_2 & \quad T^{+-}
\end{array}
$} \\  & \frame{$
\begin{array}{cc}
P_2^{+} & P_1^{+}-S_2^{+}
\end{array}
$}
\end{array}
\right\}
\end{array}
\end{equation}
This shows that to any of them the corresponding extension $\Phi _{{\cal E}
_{i-1}}$ can be applied. Returning to the sequence (\ref{r2extexts}) we see
that these extensions remove the summands of the form $S_i^{\pm }$ so that
the alternative extension also becomes applicable. As a result we get for
the Heisenberg algebra $\widetilde{\cal H}$ nine deformed costructures
$$
\begin{array}{ll}
	\Delta _{{\cal J}_1{\cal J}_0}\left( \widetilde{\cal H} \right) =
	&  \left\{
		\begin{array}{ll}
			\begin{array}{c}
			P_1^{+}\qquad  \\
			P_2^{+}\qquad
			\end{array}
		&
			\begin{array}{cc}
			T^{++} & \quad T_1\quad  \\
			T_2 & T^{++}
			\end{array}
		\\
		&
			\begin{array}{cc}
			P_2^{+} \quad  & P_1^{+} \quad
			\end{array}
		\end{array}
	\right\}
\\
\Delta _{\widetilde{{\cal E}_0}{\cal J}_1{\cal J}_0}\left( \widetilde{
\cal H}\right) = &
     \left\{
\begin{array}{ll}
\begin{array}{c}
P_1^{+} \\
P_2^{+}-S_1^{-}
\end{array}
&
\begin{array}{cc}
T^{-+} & \quad T_{1\quad } \\
T_2 & T_{R_1}
\end{array}
 \\  &
\begin{array}{cc}
P_2^{+}-S_1^{+} & P_1^{+}
\end{array}
\end{array}
\right\}
\\ \Delta _{\widetilde{{\cal E}_1}{\cal J}_1{\cal J}_0}
\left( \widetilde{\cal H}\right) = &
\left\{
\begin{array}{ll}
\begin{array}{c}
P_1^{+}-S_2^{-} \\
P_2^{+}
\end{array}
&
\begin{array}{cc}
T_{R_2} & \quad T_1\quad  \\
T_2 & T^{+-}
\end{array}
 \\  &
\begin{array}{cc}
P_2^{+} & P_1^{+}-S_2^{+}
\end{array}
\end{array}
\right\}
\\                 
\Delta _{{\cal E}_0{\cal J}_1{\cal J}_0}\left( \widetilde{\cal H}\right) =
& \left\{
\begin{array}{ll}
\begin{array}{c}
P_1^{-} \\
P_2^{+}+S_1^{-}
\end{array}
&
\begin{array}{cc}
T^{++} & \quad T_1\quad  \\
T_2 & T^{++}
\end{array}
\\  &
\begin{array}{cc}
P_2^{+}+S_1^{+} & R_1
\end{array}
\end{array}
\right\}
\\ \Delta _{\widetilde{{\cal E}_0}{\cal E}_0{\cal J}_1{\cal J}_0}\left(
\widetilde{\cal H}\right) = & \left\{
\begin{array}{ll}
\begin{array}{c}
P_1^{-}\qquad  \\
P_2^{+}\qquad
\end{array}
 &
\begin{array}{cc}
T^{-+} & \quad T_1\quad  \\
T_2 & T_{R_1}
\end{array}
 \\  &
\begin{array}{cc}
P_2^{+}\quad  & R_1\quad
\end{array}
\end{array}
\right\}
\\ \Delta _{{\cal E}_1{\cal E}_0\widetilde{{\cal E}_1}{\cal J}_1{\cal J}
_0}\left( \widetilde{\cal H}\right) = & \left\{
\begin{array}{ll}
\begin{array}{c}
P_1^{-} \\
P_2^{-}+S_1^{-}
\end{array}
&
\begin{array}{cc}
T_{R_2} & \quad T_1\quad  \\
T_2 & T^{+-}
\end{array}
\\  &
\begin{array}{cc}
R_2+S_1^{+} & R_1\;
\end{array}
\end{array}
\right\}             
\\                   
\Delta _{{\cal E}_1{\cal J}_1{\cal J}_0}\left( \widetilde{\cal H}\right) =
& \left\{
\begin{array}{ll}
\begin{array}{c}
P_1^{+}+S_2^{-} \\
P_2^{-}
\end{array}
 &
\begin{array}{cc}
T^{++} & \quad T_1\quad  \\
T_2 & T^{++}
\end{array}
 \\  &
\begin{array}{cc}
R_2 & P_1^{+}+S_2^{+}
\end{array}
\end{array}
\right\}
\\ \Delta _{{\cal E}_1{\cal E}_0\widetilde{{\cal E}_0}{\cal J}_1{\cal J}
_0}\left( \widetilde{\cal H}\right) = & \left\{
\begin{array}{ll}
\begin{array}{c}
P_1^{-}+S_2^{-} \\
P_2^{-}
\end{array}
 &
\begin{array}{cc}
T^{-+} & \quad T_1\quad  \\
T_2 & T_{R_1}
\end{array}
\\  &
\begin{array}{cc}
R_2 & R_1+S_2^{+}
\end{array}
\end{array}
\right\}
\\ \Delta _{{\cal E}_1\widetilde{{\cal E}_1}{\cal J}_1{\cal J}_0}\left(
\widetilde{\cal H}\right) = & \left\{
\begin{array}{ll}
\begin{array}{c}
P_1^{+}\qquad  \\
P_2^{-}\qquad
\end{array}
 &
\begin{array}{cc}
T_{R_2} & \quad T_1\quad  \\
T_2 & T^{+-}
\end{array}
 \\  &
\begin{array}{cc}
R_2\quad  & \;P_1^{+}\quad
\end{array}
\end{array}
\right\}
\end{array}
$$

Here is how these costructures are connected by the external and internal
extension twists:
\begin{equation}
\label{diagram}
\begin{array}{ccccc}
 \widetilde{\cal H}_{{\cal E}_1{\cal E}_0\widetilde{{\cal
E}_0}{\cal J}_1{\cal J}_0} &  &  &  &  \\
\uparrow {\Phi }_{{\cal E}_1} &  &  &  &  \\
 \widetilde{\cal H}_{{\cal E}_0\widetilde{{\cal E}_0}{\cal J}
_1{\cal J}_0} & \stackrel{\widetilde{\Phi }_{{\cal E}_0}}{\longleftarrow } &
 \widetilde{\cal H}_{{\cal E}_0{\cal J}_1{\cal J}_0} &  &
\\
\uparrow {\Phi }_{{\cal E}_0} &  & \uparrow {\Phi }_{{\cal E}_0} &  &  \\
 \widetilde{\cal H}_{\widetilde{{\cal E}_0}{\cal J}_1{\cal J}
_0} & \stackrel{\widetilde{\Phi }_{{\cal E}_0}}{\longleftarrow } &
\widetilde{\cal H}_{{\cal J}_1{\cal J}_0} & \stackrel{\widetilde{
\Phi }_{{\cal E}_1}}{\longrightarrow } &
\widetilde{\cal H}_{\widetilde{{\cal E}_1}{\cal J}_1{\cal J}_0} \\
&  & \downarrow {\Phi }_{{\cal E}_1} &  & \downarrow
{\Phi }_{{\cal E}_1}
\\  &  & \widetilde{\cal H}_{{\cal E}_1
{\cal J}_1{\cal J}_0} & \stackrel{\widetilde{\Phi }_{{\cal E}_1}}{
\longrightarrow } &
\widetilde{\cal H}_{{\cal E}_1\widetilde{{\cal E}_1}{\cal J}_1
{\cal J}_0} \\  &  &  &  & \downarrow
{\Phi }_{{\cal E}_0} \\  &  &  &  & \widetilde{\cal H}_{
{\cal E}_1{\cal E}_0\widetilde{{\cal E}_1}{\cal J}_1{\cal J}_0}
\end{array}
\end{equation}
The vertical arrows are the ordinary extension twists for Heisenberg
subalgebras. The horizontal arrows correspond to the external twists borrowed
from the chains for $U(sl(N))$. The squares of the diagram are commutative.
Its asymmetry is justified by the fact
that ${\Phi }_{{\cal E}_i}$ and  $\widetilde{\Phi }_{{\cal E}_j} $
commute only for $i=j$. At the same time in the columns of the diagram
(\ref{diagram}) both ${\Phi }_{{\cal E}_i}$ are applicable, one -- as a legal
extension for the ${\cal H}_{\cal J}$ and the other due to the corresponding
"matreshka" effect.


\section{Conclusions}
We have shown that the multidimensional Heisenberg algebras (of the type
(\ref{generators}) ) have the fixed spectrum of deformed costructures. The
properties of this spectrum is tightly connected with different chains of
extended twists that can be applied to the initial universal enveloping
algebra. It is obvious that such Heisenberg algebras can be included in
universal enveloping algebras other than $U(sl(N))$ and there they can be
treated analogously. It is equally obvious that increasing the number of
initial Jordanian twists (and enlarging correspondingly the rows and the
columns of creation and annihilation operators in the multidimensional
Heisenberg subalgebra) one can get more complicated
costructures. The main feature here is the appearance of the $S_i^{\pm}$
summands that mix the creation (annihilation) operators from different
rows (columns).


\section{Acknowledgments}

The author is heartily grateful to the organizers of the
XIV-th Max Born Symposium for their warm hospitality.

This work has been partially supported by the Russian Foundation for
Fundamental Research under the grant N 97-01-01152.


\end{document}